\documentstyle{amsppt}
\magnification=\magstep1
\NoBlackBoxes

 % triplenorm of #1

\topmatter 

\title Inequalities of correlation type for symmetric stable 
random vectors
\endtitle
\rightheadtext{Inequalities of correlation type}
\author A. L. Koldobsky and S. J. Montgomery--Smith\endauthor
\address Division of Mathematics and Statistics, 
University of Texas at San Antonio, San Antonio, TX 78249, U.S.A. 
\endaddress
\email koldobsk\@ringer.cs.utsa.edu \endemail
\address Department of Mathematics, University of Missouri--Columbia, 
Columbia, MO 65211
\endaddress
\email stephen\@mont.cs.missouri.edu \endemail

\abstract We point out a certain class of functions $f$ and $g$
for which random variables $f(X_1,\dots,X_m)$ and $g(X_{m+1},
\dots,X_k)$ are non-negatively correlated for any symmetric 
jointly stable random variables $X_i.$
We also show another result that is related to the correlation problem
for Gaussian measures of symmetric convex sets.

\endabstract

\endtopmatter \document \baselineskip=14pt                                         

\head 1. Introduction \endhead

For $0<q\le 2,$ let $Y$ be a symmetric $q$-stable
random vector in $\Bbb R^n$ with characteristic function 
$$\phi(\theta)= \exp(-\|\sum_{i=1}^n \theta_is_i\|^q),\quad 
\theta\in \Bbb R^n,
\tag{1}$$
where $s_1,\dots,s_n\in L_q([0,1]),$ and the norm is taken 
from the space $L_q([0,1]).$

For any $k\in \Bbb N$, and any choice of vectors 
$\xi_1,\dots,\xi_k\in \Bbb R^n,$ the inner products
$X_1=(Y,\xi_1),\dots, X_k=(Y,\xi_k)$ are symmetric 
$q$-stable random variables. The random variables
$X_1,\dots, X_k$ are jointly $q$-stable with zero mean,
and we say that they are $\Bbb R^n$-generated in case 
we need to emphasize the dimension of the vector $Y$.

\medbreak

In this article, we show that, for any $m<k$, and any even 
continuous positive definite
functions $f$ and $g$ on $\Bbb R^m$ and $\Bbb R^{k-m}$ respectively,
the random variables $f(X_1,\dots X_m)$ and $g(X_{m+1},\dots X_k)$
are non-negatively correlated, i.e.
$$\Bbb E\big(f(X_1,\dots, X_m)\  g(X_{m+1},\dots, X_k)\big) \ge
\Bbb E f(X_1,\dots X_m)  
\ \Bbb E g(X_{m+1},\dots, X_k),\tag{2}$$
where $\Bbb E$ stands for the expectation.

\bigbreak
Inequality (2) reminds one of some results related
to the concept of associated random variables. Recall 
that random
variables $X_1,\dots,X_k$ are said to be associated 
if, for any choice of non-decreasing (in each variable)
functions $f$ and $g$ on $\Bbb R^k,$ the random variables 
$f(X_1,\dots,X_k)$ and $g(X_1,\dots,X_k)$ are non-negatively
correlated whenever the expectations exist. 
Pitt (1982) proved that jointly Gaussian 
random variables are associated if and only if the 
correlation between each pair is non-negative.
Lee, Rachev and Samorodnitsky (1990) generalized this 
result to the case of jointly $q$-stable random variables
by giving a necessary and sufficient condition in terms of 
the spectral measure. Inequality (2) points out a special
class of functions $f$ and $g$ for which the correlation between
$f(X)$ and $g(X)$ is non-negative independently of  
relations between the jointly $q$-stable random variables $X_i$.  
For other results related to association of random variables, see
Joag-dev, Perlman and Pitt (1983), and Suquet (1994).

\medbreak

Another celebrated result of Pitt (1977) shows that, 
for any jointly Gaussian $\Bbb R^2$-generated random variables 
$X_1,\dots, X_k$, 
inequality (2) holds if $f$ and $g$ are the indicator 
functions of cubes in $\Bbb R^m$ and $\Bbb R^{k-m},$
namely, for each $t>0,$
$$P(\max_{1\le i\le k} |X_i|<t) \ge
P(\max_{1\le i\le m} |X_i|<t)        
\ P(\max_{m+1\le i\le k} |X_i|<t).\tag{3}$$
In other words, the quantity in the left-hand side
is minimal (subject to the given marginal distributions) 
if for each choice of $i,j$ with  $1\le i\le m$ 
and  $m+1\le j\le k$
the random variables $X_i$ and $X_j$ are independent,
that is to say, $b_{ij}= \hbox{\rm Cov}
(X_i,X_j)=0$.  An equivalent formulation
of the same fact is that, for any symmetric convex sets
$F$  and $G$ in $\Bbb R^2,$ $\mu(F\cap G)\ge \mu(F)\mu(G),$
where $\mu$ is a symmetric Gaussian measure in $\Bbb R^2.$
The question of whether
the same is true for symmetric convex sets in $\Bbb R^n$
( and, correspondingly, for $\Bbb R^n$-generated Gaussians)
remains open (see Schlumprecht, Schechtman and Zinn (1994)
for a historical survey and partial results).

In Section 3, we consider the quantity
in the left-hand side of (3) as a function of the $m(k-m)$
variables $b_{i,j},$ and prove that, for every dimension $n$, 
this function has a local
minimum at the origin. Note that, to solve the problem completely, 
one has to prove that the function has global minimum at the
origin.

\head 2. A correlation inequality for positive definite functions
of stable variables \endhead

In order to prove inequality (2) we need the
following simple fact.

\proclaim{Lemma 1} Let $0<q\le 2,$ and  $\xi,\eta$ be any 
vectors from the space $L_q([0,1]).$ Then
$$\exp(-\|\xi+\eta\|^q) + \exp(-\|\xi-\eta\|^q) \ge
2\exp(-\|\xi\|^q-\|\eta\|^q).$$
\endproclaim

\demo{Proof}  A result of W.~Orlicz (1933) (see also Clarkson (1936))
states that, for every $0<q\le 2$
and $\xi,\eta\in L_q,$ 
$$\|\xi+\eta\|^q+\|\xi-\eta\|^q\le 2(\|\xi\|^q+\|\eta\|^q).$$
Now use the inequality relating the arithmetic and geometric
means to obtain
$$\split \exp(-\|\xi+\eta\|^q) + \exp(-\|\xi-\eta\|^q)\ge  
\kern 3 true in\\
2\exp(-\|\xi+\eta\|^q/2-\|\xi-\eta\|^q/2)\ge 
2\exp(-\|\xi\|^q-\|\eta\|^q).\qed \endsplit $$ 
\enddemo

\medbreak

\proclaim{Theorem 1} Let $0<q\le 2$ and $X_1,\dots,X_k$ be 
jointly $q$-stable random variables. Then for any $m<k$
and any even continuous positive definite functions $f,g$ on $\Bbb R^m$
and $\Bbb R^{k-m}$\ respectively, the random variables $f(X_1,\dots,X_m)$
and $g(X_{m+1},\dots,X_k)$ are non-negatively correlated.
\endproclaim

\demo{Proof} By Bochner's theorem, $f$ and $g$ are the 
characteristic functions of finite measures
$\mu$ and $\nu$ on $\Bbb R^m$ and $\Bbb R^{k-m}$\ respectively.
The measures $\mu$ and $\nu$ are symmetric because the functions
$f$ and $g$ are even.

Let $Y$ be the $q$-stable random vector in $\Bbb R^n$
generating $X_1,\dots,X_k$, and let 
$\xi_1,\dots,\xi_k\in \Bbb R^n$ be the vectors 
for which $X_1=(Y,\xi_1),\dots, X_k=(Y,\xi_k)$.  Denote by $\gamma$
the distribution of the vector $Y,$ so $\gamma$ is a probability
$q$-stable measure in $\Bbb R^n$ with the characteristic 
function given by (1).

Using Fubini's Theorem, we see that
$$ \Bbb E\big(f(X_1,\dots, X_m)\  g(X_{m+1},\dots, X_k)\big) $$
$$ = \int_{\Bbb R^n} f\big((x,\xi_1),\dots,(x,\xi_m)\big)
\ g\big((x,\xi_{m+1}),\dots,(x,\xi_k)\big)\ d\gamma(x) $$
$$ = \int_{\Bbb R^n} \Big( \int_{\Bbb R^m} 
\exp(-i(u_1(x,\xi_1)+\dots+u_m(x,\xi_m)))\ d\mu(u_1,\dots,u_m)\ \times $$
$$ \int_{\Bbb R^{k-m}} 
\exp(-i(u_{m+1}(x,\xi_{m+1})+\dots+u_k(x,\xi_k))) 
\ d\nu(u_{m+1},\dots,u_k) \Big) \ d\gamma(x) $$ 
$$ = \int_{\Bbb R^m}\int_{\Bbb R^{k-m}} \Big(
\int_{\Bbb R^n} \exp\big(-i(x,\sum_{j=1}^k u_j\xi_j)\big)
\ d\gamma(x) \Big)
\ d\mu(u_1,\dots,u_m)\ d\nu(u_{m+1},\dots,u_k). \tag{4} $$ 

Let $\alpha= \sum_{j=1}^m u_j\xi_j,\beta=\sum_{j=m+1}^k u_j\xi_j
\in \Bbb R^n.$ Considering the coordinates of the vectors 
$\alpha$ and $\beta$ as linear functions of the coordinates of
$u_1, \dots, u_m$ and 
$u_{m+1}, \dots, u_k,$ respectively, and using (1) we see that
the quantity in (4) is equal to
$$I_1=\int_{\Bbb R^m}\int_{\Bbb R^{k-m}} 
\exp(-\|\sum_{j=1}^n \alpha_j s_j + \sum_{j=1}^n \beta_j s_j\|^q) 
\ d\mu(u_1,\dots,u_m)\ d\nu(u_{m+1},\dots,u_k), \tag{5}$$ 
where the norm is taken from the space $L_q([0,1])$.
Denote by $I_2$ the expression in (5) with minus instead of 
plus under the norm.
Since the measure $\nu$ is symmetric, $I_1=I_2.$
By Lemma 1, 
$$(I_1+I_2)/2 \ge  \int_{\Bbb R^m}\int_{\Bbb R^{k-m}} 
\exp(-\|\sum_{j=1}^n \alpha_j s_j\|^q) \ \times$$
$$\ \exp(-\| \sum_{j=1}^n \beta_j s_j\|^q) 
\ d\mu(u_1,\dots,u_m)\ d\nu(u_{m+1},\dots,u_k) $$ 
$$ = \int_{\Bbb R^m}
\exp(-\|\sum_{j=1}^n \alpha_j s_j\|^q)\ d\mu(u_1,\dots,u_m)\ \times$$ 
$$\ \int_{\Bbb R^{k-m}}  \exp(-\| \sum_{j=1}^n \beta_j s_j\|^q) 
\ d\nu(u_{m+1},\dots,u_k).$$ 
Repeating all the calculations in the reverse order
we show that the latter quantity is equal to
$\Bbb E f(X_1,\dots X_m)  
\ \Bbb E g(X_{m+1},\dots, X_k)$
which finishes the proof. \qed \enddemo

\bigbreak

\subheading{Examples} (i) Let 
$f(x_1,\dots,x_m) = (1-|x_1|)_+ \cdots (1-|x_m|)_+$, and
$g(x_{m+1},\dots,x_k) = (1-|x_{m+1}|)_+ \cdots (1-|x_k|)_+$,
where the function $(1-|t|)_{+}$ is equal to $1-|t|$ if $t\in [-1,1],$
and is equal to zero otherwise. It is well known that the function
$(1-|t|)_{+}$ is positive definite, and hence $f$\ and $g$\ are
positive definite. Thus, by Theorem 1, for every $m<k$
and every jointly stable random variables $X_1,\dots,X_k$,
$$\Bbb E\big((1-|X_1|)_{+}\cdot \dots (1-|X_k|)_{+}\big) \ge$$
$$\Bbb E\big((1-|X_1|)_{+}\cdot \dots (1-|X_m|)_{+}\big)
\  \Bbb E\big((1-|X_{m+1}|)_{+}\cdot \dots (1-|X_k|)_{+}\big).$$
The latter inequality can be generalized by taking any functions 
$f$ and $g$ of the form  $f(x_1,\dots,x_m)=f_1(x_1)\dots f_m(x_m),
\ g(x_{m+1},\dots,x_k)=f_{m+1}(x_{m+1})\dots f_k(x_k),$ where
$f_1,\dots,f_k$ are even functions on $\Bbb R$ which are convex 
and decreasing on $[0,\infty).$ Such functions $f_i$ are positive
definite by a well-known result of Polya.
\medbreak
(ii)\ Let $q_1,\dots,q_k\in (0,2]$, 
$f(x_1,\dots,x_m) = \exp(-|x_1|^{q_1} - \cdots - |x_m|^{q_m})$, and
$g(x_{m+1},\dots,x_k) = \exp(-|x_{m+1}|^{q_{m+1}} - \cdots - |x_k|^{q_k})$.
Since for any $q \in (0,2]$\ the function $\exp(-|t|^q)$\ is positive
definite, it follows that $f$\ and $g$\ are positive definite.
Therefore, for every $m<k$
$$\Bbb E\big(\exp(-|X_1|^{q_1}-\dots-|X_k|^{q_k})\big) \ge$$
$$\Bbb E\big(\exp(-|X_1|^{q_1}-\dots-|X_m|^{q_m})\big)
\Bbb E\big(\exp(-|X_{m+1}|^{q_{m+1}}-\dots-|X_k|^{q_k})\big).$$

\medbreak

\subheading{Remarks} 
(i) In the case of jointly Gaussian random variables
the result of Theorem 1  can be extended
to some classes of continuous  functions $f$ and $g$ with power
growth at infinity and such that their Fourier transforms 
(in the sense of distributions) are non-negative locally integrable
functions with power growth at infinity. 
To do that, consider the convolutions
of the functions $f$ and $g$ with Gaussian densities $e_n$
approaching the $\delta$-function as $n\to \infty,$ and
slightly modify the proof of Theorem 1. 
\smallbreak
(ii) Y. Hu has recently 
proved that, for any even convex functions $f$ and $g$ on $\Bbb R^n$
and jointly Gaussian random variables $X_1,\dots,X_n,$ the 
random variables $f(X_1,\dots,X_n)$ and $g(X_1,\dots,X_n)$
are non-negatively correlated (private communication from 
T. Schlumprecht; compare the result of Hu with our Example 1).

\bigbreak

\head 3. On the local minimum in the correlation
for Gaussian measures of symmetric convex sets   \endhead

Let $\nu$ be the standard symmetric Gaussian measure on
$\Bbb R^n$.  Is it true that
$$\nu(F\cap G)\ge \nu(F)\nu(G)\tag{6}$$ 
for all symmetric convex sets $F$ and $G$ in $\Bbb R^n$~?
In 1977, L. Pitt proved that the answer is positive in the case 
$n=2$.  However, the question of whether the answer is
positive for every dimension $n$ is still open.

It can be seen that it suffices to consider the sets
$F=\{x\in \Bbb R^n:\ |(x,\xi_1)|\le 1,\dots,|(x,\xi_k)|\le 1\}$ and
$G=\{x\in \Bbb R^n:\ |(x,\xi_{k+1})|\le 1,\dots,|(x,\xi_{2k})|\le 1\}$,
where $k$ is an integer, and
$\xi_1,\dots,\xi_k$, $\xi_{k+1},\dots,\xi_{2k}\in \Bbb R^n$.
For these sets $F$ and $G,$ inequality (6) can be written in 
the form
$$P(\max_{1\le i\le 2k} |X_i|<1) \ge
P(\max_{1\le i\le k} |X_i|<1)        
\ P(\max_{k+1\le i\le 2k} |X_i|<1), \tag{7}$$
where $X_1,\dots, X_{2k}$ are the jointly Gaussian
random variables generated by the vectors 
$\xi_1,\dots,\xi_k$, $\xi_{k+1},\dots,\xi_{2k}\in \Bbb R^n$
and a standard Gaussian random vector $Y$ in $\Bbb R^n$, 
so that $X_i=(Y,\xi_i)$
for each $i.$

It is easy to see that, to prove inequality (6),
it suffices to consider the case where the vectors $\xi_i,\ i=1,\dots,2k$
are linearly independent. For example, if $n<2k$ and
the system of vectors $\xi_i$ has rank $n$, 
we can transfer everything to the
space $\Bbb R^{2k},$ and  
consider the vectors $\eta_i=\xi_i + \epsilon e_i
\in \Bbb R^{2k},\ i=1,\dots,2k $ 
where, for each $i,$ 
either $e_i=0$ or $\|e_i\|=1$ and $e_i$ is orthogonal to
each of the vectors $\xi_j,\ j=1,\dots,2k$ and $e_j,\ j\neq i,$
\ so that the vectors $\eta_i$ 
are linearly independent in $\Bbb R^{2k}.$ 
Then inequality (7) for the random variables generated by
the vectors $\eta_i$ would
imply inequality (7) for the random variables 
generated by $\xi_i$'s by taking the limit as 
$\epsilon \to 0$ and applying the Lebesgue dominated convergence 
theorem.

Assume that the vectors $\xi_i\in \Bbb R^{2k},\ i=1,\dots,2k$ are
linearly independent. Then the joint distribution $\mu$
of random variables $X_1,\dots,X_{2k}$ is a non-singular
Gaussian measure in $\Bbb R^{2k},$ and the left-hand side
of (7) is equal to
$$P(\max_{1\le i\le 2k} |X_i|<1) = \mu([-1,1]^{2k}).$$

We fix the scalar products $(\xi_i,\xi_j)$ for all choices of $i,j$ 
with  either $1\le i,j\le k$ or $k+1\le i,j \le 2k,$\ and
consider the quantity  $\mu([-1,1]^{2k})$ as a function of $k^2$
variables $b_{i,j}=\hbox{\rm Cov}
(X_i,X_j)$, $i=1,\dots,k$, $j=k+1,\dots,2k$.
To prove Pitt's inequality, one has to show that this
function has a global minimum at zero. Being unable to do that
we show instead that the function has a local minimum 
at zero.  This fact is a simple consequence of Theorem 2 below.

\medbreak

\def\half{\textstyle{1\over2}}
\def\d{\partial}
In the proof of Theorem 2 we use one result about log-concave functions.
A non-negative function $f$ on $\Bbb R^k$ is called log-concave if, for
every choice of $x,y\in \Bbb R^k,$ and $0\le t\le 1,$
$$f(tx+(1-t)y)\ge f(x)^t f(y)^{1-t}.$$
This means that the function $\log(f)$ is concave. 
Prekopa (1973) and Leindler (1972) have proved that
if $f$ is a log-concave function on $\Bbb R^k$ and
$0<m<k,$ then the function
$$g(x_1,\dots,x_m)=\int_{\Bbb R^{k-m}} f(x_1,\dots,x_m,z_1,\dots,
z_{k-m})\ dz$$
is also log-concave.

\proclaim{Theorem 2} Let $F$ and $G$ be symmetric convex
sets in $\Bbb R^k,$ and $\mu_B$ be a non-singular probability
Gaussian measure
in $\Bbb R^{2k}$ with the covariance matrix 
$\Cal A= \left[\matrix A & B \cr B^T & C\cr  \endmatrix \right].$
Fix the $k\times k$ matrices $A$ and $C,$ and 
consider $B=(b_{i,j})_{i,j=1}^k$ as a variable from 
the space $\Bbb R^{k^2}.$ Then the function 
$B\mapsto \mu_{B}(F \times G)$ has a local minimum at 
the point $B=0.$
\endproclaim

\demo{Proof} 
Without loss of generality, we may suppose that $F$\ and $G$\ have compact
closure.
Let $\chi_F,\ \chi_G$ be the indicator functions 
of the sets $F$ and $G.$ Taking Fourier transforms, we obtain
$$ \mu_{B}(F\times G) =
\int_{\Bbb R^k} \int_{\Bbb R^k} \chi_F(x) \chi_G(y)\, d\mu_B(x,y)$$
$$ = \int_{\Bbb R^k} \int_{\Bbb R^k} \hat{\chi_F}(x) \hat{\chi_G}(y)
   \exp(-\half(x^T A x + y^T C y + 2 x^T B y)) \, dx \, dy .$$
Taking the second partial derivative by $b_{i,j}$ and $b_{m,n},$ we get
$$ H_{i,j,m,n} =
   {\d^2\over \d b_{i,j} \d b_{m,n}} \mu_{B}(F\times G)$$
$$=\int_{\Bbb R^k} \int_{\Bbb R^k} \hat{\chi_F}(x) \hat{\chi_G}(y)
(x_i x_m y_j y_n) \exp(-\half(x^T A x + y^T C y + 2 x^T B y)) \, dx \, dy$$
$$={1\over (2 \pi)^{k} |\Cal A|^{1/2}}
\int_F \int_G {\d^4\over \d x_i \d x_m \d y_j \d y_n}
\exp(-\half (x,y)^T \Cal A^{-1} (x,y)) \, dy \, dx .$$
The fact that $|\Cal A|\neq 0$, and the validity of using 
Parseval's Equality in the latter equations,
follow from the non-singularity of the measure $\mu_B.$

Since the sets $F$ and $G$ are symmetric, the partial derivative 
of the function $B\mapsto \mu_{B}(F \times G)$ by each $b_{i,j}$
is equal to zero at the point $B=0.$
In order to show that there is a local minimum at $B=0$, we need to know
that $H$\ is positive definite when $B = 0$.  Furthermore, by a change
of variables, we see that it is sufficient to consider the special
case when $A = C = I$.  Hence, we need to show
the positive definiteness of
$$ H_{i,j,m,n} = {1\over (2 \pi)^{2k}} L_{i,m} K_{j,n} $$
where 
$$ L_{i,m} = \int_F (x_i x_m - \delta_{i,m}) \exp(-\half x^Tx) \, dx $$
and
$$ K_{j,n} = \int_G (y_j y_n - \delta_{j,n}) \exp(-\half y^Ty) \, dy . $$
Since $H = L \otimes K$, it is sufficient to show that $L$\ and $K$\
are negative definite, and clearly it is enough just to prove it for
$L$.

Thus we desire to show that 
$$ \sum_{i,m} L_{i,m} \alpha_i \alpha_m
   =
   \int_F ((\sum_i \alpha_i x_i)^2 - \|\alpha\|_2^2) \exp(-\half x^Tx) \, dx
   < 0 $$
for all $\alpha \ne 0$.  But by a change of variables, it is sufficient
to show 
$$ \int_F (x_1^2 - 1) \exp(-\half x^Tx) \, dx < 0 $$
for every convex symmetric set $F$\ with compact closure. 

To show this, we see this as
$$ \int_{-\infty}^\infty (x_1^2 - 1) \exp(-\half x_1^2) \phi(x_1) \, dx_1 ,$$
where
$$ \phi(x_1) = \int_{\Bbb R^{k-1}} \chi_F(x_1,\dots,x_k) 
\exp(-\half(x_2^2+\dots+x_k^2)) \, dx_2\dots dx_k .$$
Since $\chi_F(x) \exp(-\half(x_2^2+\dots+x_k^2))$\ is log-concave
in $\Bbb R^k,$ the
result of Prekopa and Leindler mentioned before the formulation
of Theorem 2 implies that $\phi$\ is also log-concave.  Since
$\phi$\ is also symmetric, it follows that
$\phi(x_1) = \phi_1(|x|)$, where $\phi_1$\ is a decreasing function. 
Furthermore, since $F$\ has compact closure, $\phi_1$\ is non-constant.
Hence in order to show that
$$ \int_{-\infty}^\infty (x_1^2 - 1) \exp(-\half x_1^2) \phi(x_1) \, dx_1 < 0,$$
it is sufficient to show that for all $0<a<\infty$
$$ \theta(a) = \int_{-a}^a (x_1^2 - 1) \exp(-\half x_1^2) \, dx_1 < 0 .$$
The function under the latter integral has antiderivative
$-x_1 \exp(- \half x_1^2),$ so the result follows.\qed \enddemo

\bigbreak

Finally, we present one more argument showing that 
inequality (6) would be proved if one showed that
the function from Theorem 2 had 
global minimum at zero. 
\smallbreak
Let $A = C = I$. Since the sets $F$\ and $G$\ are convex, 
their topological boundaries have
zero Lebesgue measure. Let $\nu$\ be standard Gaussian measure on
$\Bbb R^k$.  Then $\mu_0(F \times G) = \nu(F) \nu(G)$, whereas
$\lim_{\lambda \to 1} \mu_{\lambda I}(F \times G) = \nu(F \cap G)$.
To see this last assertion, note that
$$ \eqalignno{
   \mu_{\lambda I}(F \times G)
   &= {1 \over ((2 \pi(1-\lambda^2))^k}
     \int_F \int_G \exp(-{1\over 2(1-\lambda^2)} 
                   (x^T x - 2\lambda x^T y + y^T y)) \, dy \, dx \cr
\noalign{\noindent which, making the substitution $x = u+v$, $y=u-v$}
   &= {1\over (\pi (1-\lambda^2))^k}
    \int_{\Bbb R^k} \int_{(F-u)\cap(u-G)}
    \exp(-{u^2\over 1+\lambda} - {v^2\over 1-\lambda}) \, dv \, du .\cr}$$
Now, if $u$\ is not in the boundary of $F$\ or the boundary of $G$, then
it is easily seen that
$$ \lim_{\lambda\to 1}
   {1\over (\sqrt{\pi} (1-\lambda))^k}
   \int_{(F-u)\cap(u-G)} \exp(-{v^2\over 1-\lambda}) \, dv
   = \chi_{F \cap G}(u) .$$
Hence the last assertion follows by Lebesgue's law of dominated convergence.

It is clear now that, if the function $\mu_B$ has global minimum
at zero then $\mu_{\lambda I}(F \times G)\ge \mu_0(F \times G),$
and, hence, $\nu(F\cap G)\ge \nu(F)\nu(G).$ 
However, the question of whether the function from Theorem 2
has global minimum at zero remains open.
\bigbreak

\subheading{Acknowledgements} We would like to thank
T. Schlumprecht, G. Schechtman and J. Zinn for bringing
the problem to our attention and providing us with
updated information including their unpublished results.

\enddocument